\documentclass[11pt]{amsart}

\usepackage{enumerate,url,amssymb,  mathrsfs, graphicx, pdfsync}
\newtheorem{theorem}{Theorem}[section]

\newtheorem*{lemma*}{Lemma}

\newtheorem{proposition}[theorem]{Proposition}
\newtheorem{corollary}[theorem]{Corollary}
\theoremstyle{definition}

\newtheorem{example}[theorem]{Example}

\theoremstyle{remark}

\numberwithin{equation}{section}

\renewcommand{\bar}[1]{\overline{#1}}

\newcommand{\dd}{\mathbb{D}}
\newcommand{\rr}{\mathbb{R}}

\newcommand{\abs}[1]{\lvert#1\rvert}

\newcommand{\C}{\mathbb{C}}

\newcommand{\W}{\mathscr{W}}

\renewcommand{\P}{\mathcal{P}}

\newcommand{\R}{\mathbb{R}}
\newcommand{\X}{\mathbb{X}}

\newcommand{\Y}{\mathbb{Y}}

\newcommand{\dtext}{\textnormal d}

\newcommand{\onto}{\xrightarrow[]{{}_{\!\!\textnormal{onto\,\,}\!\!}}}

\DeclareMathOperator{\loc}{loc}

\def\leq{\leqslant}
\def\geq{\geqslant}

\def\le{\leqslant}
\def\ge{\geqslant}

\def\XXint#1#2#3{{\setbox0=\hbox{$#1{#2#3}{\int}$}\vcenter{\hbox{$#2#3$}}\kern-.5\wd0}}

\def\XXiint#1#2#3{{\setbox0=\hbox{$#1{#2#3}{\iint}$}\vcenter{\hbox{$#2#3$}}\kern-.5\wd0}}

\begin{document}
\title{Bi-Sobolev extensions}

\author[A. Koski]{Aleksis Koski}
\address{Departamento de Matem\'aticas, Univeridad Aut\'onoma de Madrid, E-28049 Madrid, Spain}
\email{aleksis.koski@gmail.com}

\author[J. Onninen]{Jani Onninen}
\address{Department of Mathematics, Syracuse University, Syracuse,
NY 13244, USA and  Department of Mathematics and Statistics, P.O.Box 35 (MaD) FI-40014 University of Jyv\"askyl\"a, Finland
}
\email{jkonnine@syr.edu}

\thanks{A. Koski was supported by by the ERC Advanced Grant number 834728.
J. Onninen was supported by the NSF grant  DMS-2154943.   
}

\subjclass[2010]{Primary 46E35, 30C62. Secondary 58E20}


\keywords{Sobolev homeomorphisms, Sobolev extensions, Beurling-Ahlfors extension, harmonic extension, quasiconformal mapping and mapping of finite distortion.}

\begin{abstract} 
We give a full  characterization of  
circle homeomorphisms which  admit a homeomorphic  extension to the unit disk  with finite bi-Sobolev norm.
As a special case,  a  bi-conformal  variant  of the famous Beurling-Ahlfors extension theorem is obtained. Furthermore we show that  the existing extension techniques  such as applying either the harmonic  or  the Beurling-Ahlfors  operator work poorly in the degenerated setting.  This also gives an affirmative  answer to a  question of Karafyllia and Ntalampekos.
\end{abstract}

\maketitle
\section{Introduction}
The celebrated Beurling-Ahlfors quasiconformal extension theorem~\cite{AhBe} states that a self-homeomorphism of the unit disk $\dd \subset \mathbb C$   is quasiconformal if and only if the boundary correspondence  mapping $\varphi \colon \mathbb S \onto \mathbb S$ is quasisymmetric. It has found a number of  applications in Teichm\"uller theory, Kleinian groups, conformal welding
and dynamics, see e.g.~\cite{AIMb, Hu}.  Recall that a sense-preserving homeomorphism $\varphi \colon \mathbb S \onto \mathbb S$ is {\it quasisymmetric} if the quantity
\begin{equation}\label{eq:pointquasi}
\delta_\varphi (\theta, t) = \max \left\{  \frac{\varphi (e^{i (\theta +t)})  - \varphi (e^{i\theta })   }{\varphi (e^{i\theta }) - \varphi (e^{i (\theta -t)}) }  \, , \,   \frac{  \varphi (e^{i\theta }) - \varphi (e^{i (\theta -t)})   }{        \varphi (e^{i (\theta +t)})  - \varphi (e^{i\theta })        }  \right\}
\end{equation}
is uniformly bounded both in $\theta \in [0, 2 \pi]$ and in $t \in (0, 2 \pi)$. Furthermore,  a  homeomorphism $h \colon \mathbb D \onto \mathbb D$ of Sobolev class $\W_{\loc}^{1,1} (\mathbb D  , \mathbb C)$ is {\it quasiconformal} if 
\begin{equation}\label{eq:distortio}
\abs{Dh(x)}^2 \le K(x) J_h (x) 
\end{equation}
for some  $K\in \mathscr L^\infty (\mathbb D)$.  Here, $\abs{Dh(x)}$ is the operator norm of the weak differential $Dh(x) \colon \mathbb D \to \R^2$ of $h$ at a point $x \in \mathbb D$, and $J_h(x)= \det Dh(x)$ is  the \textit{Jacobian determinant} of $h$.  The smallest function $K(x)=K_h(x) \ge 1$ for which the distortion inequality~\eqref{eq:distortio} holds is called the {\it distortion} of $h$. 

 The remarkable feature  of a quasiconformal mapping is that its inverse is also quasiconformal. In particular, both the mapping $h \colon \dd \onto \dd$ and its inverse $f=h^{-1} \colon \dd \onto \dd$ have finite {\it conformal energy} (also called the {\it Dirichlet energy}). Their sum
\begin{equation}
\mathsf E [h] = \int_\dd \abs{Dh(x)}^2\, \dtext x + \int_\dd \abs{Df(y)}^2 \, \dtext y < \infty 
\end{equation}
is called {\it bi-conformal energy} of $h$ (and $f$). An orientation-preserving  homeomorphism $h\colon \dd \onto \dd$  with finite bi-conformal energy is  a {\it mapping of bi-conformal energy}.  Note that a homeomorphism $h\colon \mathbb D \onto \mathbb D$ is a mapping of bi-conformal energy  if and only if  the mapping $h$ lies in  $\W^{1,2} (\mathbb D, \C)$ and  has integrable distortion; that is, $K_h\in \mathscr L^1 (\mathbb D)$, see~\cite{AIMO, HKb, HKO, IOZ}. For such mappings we have
\begin{equation}\label{eq:biconf}
\begin{split}
\mathsf E [h] & = \int_\dd \abs{Dh(x)}^2\, \dtext x + \int_\dd \abs{Df(y)}^2 \, \dtext y \\&=   \int_\dd K_f(y)\, \dtext y + \int_\dd K_h(x) \, \dtext x \\ &= \int_\dd \left[  \abs{Dh (x)}^2 + K_h(x)  \right] \, \dtext x < \infty \, . 
\end{split}
\end{equation}
The last integral is also well defined  in the class of non-injective mappings of
integrable distortion.  Such deformations and their $n$-dimensional counterparts play an important role in   Geometric Function Theory (GFT)~\cite{AIMb, HKb, IMb} as they share several fundamental topological and analytical
properties of analytic functions.  For instance, a nonconstant mapping $h$ with   $\mathsf E [h] <\infty$ is continuous, discrete and open~\cite{IS}.   As a result,  mappings of bi-conformal energy form the widest class of homeomorphisms that one can hope to build a viable extension of GFT with connections to mathematical models of Nonlinear Elasticity~\cite{Anb, Bac, Cib}.    This circle of ideas has applications in materials science and critical phase phenomena, where the distortion functionals are natural measures of change in a system.  They address fundamental questions of microstructure and length scales. 
In the setting of unbounded distortion, however, many of the usual tools of quasiconformal mappings are lost. In particular, this degenerated setting truly challenges the available extension methods in the literature. Nevertheless,  it is  known that the Beurling-Ahlfors  extension operator can be pushed to  produce
a homeomorphism $h \colon \dd \onto \dd$ with finite    bi-conformal energy~\eqref{eq:biconf}  if  the boundary  correspondence   homeomorphism $\varphi \colon \mathbb S  \onto \mathbb S $ satisfies  the condition
\begin{equation}\label{eq:sufficient}
\int_{\mathbb S} \int_{\mathbb S}   \left| \frac{\varphi (\xi) - \varphi (\eta)}{ \xi - \eta }\right|^2   {\dtext \xi } \,  {\dtext \eta }  +  \int_0^{2\pi }  \int_0^{2\pi }  \delta_\varphi (\theta, t)\,  \dtext \theta \,  \dtext t   < \infty \, . 
\end{equation}
This follows from a very recent result  of  Karafyllia and Ntalampekos~\cite{KN}. Actually, their sharp $\mathscr L^\infty$-estimates for  the Beurling-Ahlfors   operator also give sufficient conditions for 
 other extension problems involving unbounded distortion functions,  see also~\cite{CCH, QZ, Sa, Za}. Such extension problems are in demand e.g. in complex dynamics,~\cite{Ha, PZ}.  In particularly,  Karafyllia and Ntalampekos showed that  a boundary homeomorphism $\varphi \colon \mathbb S \onto \mathbb S$ admits a homeomorphic extension $h \colon \mathbb D \onto \mathbb D$ with  $p$-integrable distortion; that is, $K_h \in \mathscr L^p (\mathbb D)$, $p\ge 1$, if
\begin{equation}\label{eq:p-integralbdd}
\int_0^{2\pi }  \int_0^{2\pi }  \left[ \delta_\varphi (\theta, t)\right]^p \,  \dtext \theta \,  \dtext t   < \infty \, . 
\end{equation}  
 They raised a question  if the sufficient condition~\eqref{eq:p-integralbdd} is also necessary for obtaining an extension of $p$-integrable distortion, see \cite[Question 1.5]{KN}. Our next result, however, shows that this is far from being the case.
 \begin{example}\label{ex:p-integrabledist}
For $p\ge 1$ and every $q>p$,  there exists a Lipschitz homeomorphism $h \colon \overline{\mathbb D} \onto \overline{\mathbb D}$ such that $K_h \in \mathscr L^p (\mathbb D)$ and 
\[ \int_0^{2\pi }  \int_0^{2\pi } \log ^{q} \big( e+  \delta_\varphi (\theta, t) \big) \,  \dtext \theta \,  \dtext t   = \infty \, ,   \]
  where $\varphi =h$ on $\partial \mathbb D$.
 \end{example}
 In particular, this provides  us a mapping of bi-conformal energy from the unit disk onto itself whose boundary mapping does not satisfy the  condition~\eqref{eq:sufficient}. Nevertheless,  the first part of~\eqref{eq:sufficient} still remains  as a necessary requirement.  Indeed, a mapping $\varphi \colon \mathbb S \to \mathbb C$ admits a continuous extension (not necessarily  a homeomorphism) to the unit disk $\mathbb D$ in the Sobolev class $\W^{1,p} (\mathbb D , \mathbb C)$, $1<p<\infty$ if and only if   it satisfies the so-called {\it $p$-Douglas condition},
\begin{equation}\label{eq:pdouglas}
\int_{\mathbb S} \int_{\mathbb S} \left| \frac{\varphi (\xi) - \varphi (\eta)}{ \xi - \eta }\right|^p {\dtext \xi } \,  {\dtext \eta }  < \infty \, ,  
\end{equation}
known as the {\it Douglas condition} when p=2,~\cite{Do}. For the proof of this result we refer to~\cite[p. 151-152]{St}.  Note that an arbitrary boundary homeomorphism $\varphi \colon \mathbb S \onto \mathbb S$ satisfies the $p$-Douglas condition for $p<2$. 

On the other hand,  according to the famous theory of Rad\'o~\cite{Rad}, Kneser~\cite{Kn} and Choquet~\cite{Ch}, see also~\cite{Dub},  any homeomorphic boundary value $\varphi \colon \mathbb S \onto \mathbb S $ admits a homeomorphic harmonic extension of $\mathbb D$ onto itself. Moreover,
since  the harmonic extension of $\varphi$  has the smallest Dirichlet energy among all extensions, a homeomorphism $\varphi \colon \mathbb S \onto \mathbb S$ admits a homeomorphic extension $h \colon \overline{\dd} \onto \overline{\dd}$ in $\W^{1,2} (\dd, \C)$ if and only if $\varphi$ satisfies the Douglas condition.
Similarly for $p>1$, the $p$-harmonic variants of the Rad\'o-Kneser-Choquet theorem~\cite{AS}  show that if  a homeomorphism $\varphi \colon \mathbb S \onto \mathbb S$  satisfies the $p$-Douglas condition, then it actually admits  a homeomorphic extension $h \colon \overline{\dd} \onto \overline{\dd}$ in $\W^{1,p} (\dd, \C)$.  For $p<2$, this can be also seen  by simply  applying the standard harmonic extension operator. Indeed,  the harmonic extension of an arbitrary  homeomorphism $\varphi \colon \mathbb S \onto \mathbb S$ lies in  $\W^{1,p} (\mathbb D, \C)$ for all $p\in [1,2)$, see~\cite{Ve}. The harmonic extension operator  is actually surprisingly robust. It also works above its natural $\W^{1,2}$-domain of definition.
\begin{theorem}\label{thm:rkcp>2}
Let $\varphi \colon \mathbb S \onto \mathbb S$ be a homeomorphism and $h \colon \overline{\mathbb D } \onto \overline{\mathbb D }$ the harmonic   extension of $\varphi$.  Then the homeomorphism $h$  lies in $\W^{1,p} (\mathbb D , \mathbb C)$ if and only if the boundary map $\varphi$ satisfies the $p$-Douglas condition with $p\in (1, \infty)$.
\end{theorem}

In particular, a necessary condition for a homeomorphism $\varphi \colon \mathbb S \onto \mathbb S$ to be the trace of mapping of bi-conformal energy is that both $\varphi$ and $\varphi^{-1}$ satisfy the Douglas condition. In that case we say that $\varphi$ enjoys   the {\it bi-Douglas condition},
\begin{equation}\label{eq:bidouglas}
\int_{\mathbb S} \int_{\mathbb S}  \left| \frac{\varphi (\xi) - \varphi (\eta)}{ \xi - \eta }\right|^2   +   \left| \frac{\varphi^{-1} (\xi) - \varphi^{-1} (\eta)}{ \xi - \eta }\right|^2  {\dtext \xi } \,  {\dtext \eta } < \infty \, . 
\end{equation}
The bi-Douglas condition can be written purely in terms of the mapping $\varphi$, see~\cite{AIMO}. Indeed, a homeomorphism $\varphi \colon \mathbb S \onto \mathbb S$ satisfies the bi-Douglas condition if and only if
\begin{equation}\label{eq:bidouglas2}
\int_{\mathbb S} \int_{\mathbb S}     \left| \frac{\varphi (\xi) - \varphi (\eta)}{ \xi - \eta }\right|^2 +   \big| \log \abs{\varphi (\xi) - \varphi (\eta)    } \big|      {\dtext \xi } \, {\dtext \eta } < \infty \, . 
\end{equation}
In spite of that the harmonic extension operator is a powerful tool to produce Sobolev homeomorphic extensions (Theorem~\ref{thm:rkcp>2}),
 this techniques  certainly has its limitations  when one, for instance, moves to the bi-Sobolev setting.
 \begin{example}\label{ex:biharm}
 There exists a  homeomorphism $\varphi \colon \mathbb S \onto \mathbb S$ which satisfies the bi-Douglas condition~\eqref{eq:bidouglas} such that the inverse of the harmonic extension of $\varphi$ does not belong to $\W^{1,2} (\dd, \C)$. 
 \end{example}
Thus,  Example~\ref{ex:p-integrabledist} and Example~\ref{ex:biharm} show that new and direct ways of constructing  homeomorphic extensions are needed to prove  the following bi-conformal variant of  the Beurling-Ahlfors extension theorem. 
\begin{theorem}\label{thm:p=2=q}
Let $\varphi \colon \mathbb S \onto \mathbb S$ be an orientation-preserving homeomorphism. Then $\varphi$ satisfies the bi-Douglas condition if and only if it admits a bi-conformal energy extension to $\dd$.
\end{theorem}
This actually follows as a special case ($p=2=q$) from our main result which completely characterizes circle homeomorphisms that admit bi-Sobolev extensions to the unit disk.
\begin{theorem}\label{thm:main}
Let $\varphi \colon \mathbb S \onto \mathbb S$ be a homeomorphism. Then $\varphi$ satisfies the $p$-Douglas condition and $\varphi^{-1}$  the $q$-Douglas condition if and only if the mapping $\varphi$ admits a homeomorphic extension $h \colon \overline{\mathbb D} \onto \overline{\mathbb D}$ so that $h \in \W^{1,p} (\mathbb D , \mathbb C)$ and $h^{-1} \in \W^{1,q} (\mathbb D , \mathbb C)$.
\end{theorem}
Since every homeomorphism $\varphi  \colon \mathbb S \onto \mathbb S$ satisfies the $p$-Douglas condition with $p<2$, the following result is another immediate consequence of Theorem~\ref{thm:main}.
\begin{corollary}
Any homeomorphism $\varphi \colon \mathbb S \onto \mathbb S$ admits a homeomorphic extension to $\mathbb D$ so that both the mapping and its inverse lie in the Sobolev class $\W^{1,p} (\mathbb D, \mathbb C)$ for all $p\in [1,2)$.
\end{corollary}

The {\it principle of non-interpenetration of matter} in the mathematical models of Non-linear Elasticity   asserts that the energy-minimal displacement field $h \colon \X \onto \Y$ is a homeomorphism.  Clearly, the validity of this principle  depends on the studied {\it stored energy} functional 
\begin{equation}\label{energ}
\mathcal E_\mathbb X [h] =  \int_\mathbb X \mathbf {\bf E}(x,h, Dh )\,  \dtext x\, ,   \qquad  \mathbf E \colon  \mathbb X \times \mathbb Y \times \mathbb R^{2 \times 2} 
\end{equation}
where   the so-called \textit{stored energy function} $\mathbf E$ characterizes the mechanical and elastic properties of the material occupying the domains $\X$ and $\Y$.
The \textit{$p$-harmonic energy}, including the Dirichlet integral ($p=2$),  and the so-called  \textit{total $p$-harmonic energy} serve  as model examples,
\[
\mathscr E_p[h] = \int_\mathbb X \abs{Dh(x)}^p \, \dtext x \quad  \textnormal{and} \quad \mathscr T_p [h] =  \int_\mathbb X |Dh(x)|^p\,  \dtext x  + \int_\mathbb Y |Dh^{-1}(y)|^p  \, \dtext y\, . 
\]
For instance,  in the  case of  the Diriclet energy  the injectivity  may be lost when passing to the limit of the minimizing sequence of homeomorphisms, the {\it interpenetration of matter} takes place~\cite{IOmhh}.   Of course, here we  need to know that the minimization problem is well-posed to start with. This, and more generally  minimizing the $p$-harmonic energy,   leads  us to study the associated  Sobolev variants of the Jordan-Sch\"onflies theorem.
The classical {\it Jordan-Sch\"onflies theorem} states that for a given pair of  (simply connected) Jordan domains $\X, \Y \subset \mathbb \R^2$ and a boundary homeomorphism $\varphi \colon  \partial \mathbb X \onto \partial \Y$ there is a homeomorphism $h \colon \overline{\mathbb X }\to  \overline{\Y}$ which equals $\varphi$ on $\partial \X$.

For $p\ge 1$, the  {\it Sobolev Jordan-Sch\"onflies problem} asks to characterize the pairs of Jordan domains $\X, \Y \subset \R^2 $  for which  any boundary  homeomorphism $\varphi \colon \partial \mathbb X \onto \partial \mathbb Y$ that admits a continuous extension to $\overline{\mathbb X }$ of the Sobolev class $\W^{1,p}(\mathbb X, \R^2)$ also admits a homeomorphic extension $h \colon \mathbb X \onto \Y$  in  $\W^{1,p} (\mathbb X, \R^2)$. The Sobolev Jordan-Sch\"onflies problem is well understood today~\cite{KKO, KOext1, KOext3}. On the other hand,  the well-posedness of variational problems involving the total $p$-harmonic energy requires us to study the associated bi-Sobolev Jordan-Sch\"onflies  problems.

Since the Sobolev spaces under consideration are invariant via a global bilipschitz
change of variables in both the reference and the deformed configuration, Theorem~\ref{thm:main} generalizes right away
to the case of Lipschitz domains (an axiomatic assumption in the theory of Non-linear Elasticity). Thus, the pair of Lipschitz domains satisfies the conditions in the bi-Sobolev Jordan-Sch\"onflies  problem. Moreover, a standard reflection argument shows that the extension given by Theorem \ref{thm:main} may also be defined globally with $h \in W^{1,p}_{\loc}(\mathbb C, \mathbb C)$ and $h^{-1} \in W^{1,q}_{\loc}(\mathbb C, \mathbb C)$. Thus we obtain the following theorem as a corollary.
\begin{theorem}\label{thm:lipvariant}
Let $\X$ and $\Y$ be simply connected bounded Lipschitz domains in the complex plane. Suppose that a homeomorphism $\varphi \colon \partial \X \onto \partial \Y$ admits a continuous extension  $f \colon \overline{\X} \to \mathbb C$ which lies in the Sobolev space $\W^{1,p} (\X, \mathbb C)$ for some $1\le p < \infty$. If $\varphi^{-1}$ also admits a continuous extension to $\Y$ in  $\W^{1,q} (\Y, \mathbb C)$ for some $1\le q < \infty$, then $\varphi$ admits a homeomorphic extension $h \colon \mathbb C \onto \mathbb C$ in $\W_{\loc}^{1,p} (\mathbb C, \mathbb C)$ so that the inverse $h^{-1}\in \W_{\loc}^{1,q} (\mathbb C, \mathbb C)$. 
\end{theorem}
Now, if an orientation-preserving  boundary homeomorphism $\varphi \colon \partial \X \onto \partial \Y$ satisfies the assumptions of Theorem~\ref{thm:lipvariant} with $p=2=q$, then a minimizer of   the bi-conformal energy
\[
\mathsf E [h]  = \int_\X \abs{Dh(x)}^2\, \dtext x + \int_\Y \abs{Dh^{-1}(y)}^2 \, \dtext y = \int_\X \left[  \abs{Dh (x)}^2 + K_h(x)  \right] \, \dtext x  
\]
exists among all homeomorphisms $h \colon \overline{\X} \onto \overline{\Y}$ which coincide with $\varphi$ on the boundary of $\X$; that is, an interpenetration of matter does not occur. It is worth comparing this to  the mappings with smallest Dirichlet energy which need not be injective  in the case of   nonconvex target domains~\cite{AN, Ch, IOinnerv}.

\section{Preliminaries}

For the proofs to be presented in the next sections we will need an alternative formulation of the $p$-Douglas condition \eqref{eq:pdouglas}. To state this condition, we let $I_{n,k}$ denote a dyadic decomposition of the unit circle. Then the result reads as follows.

\begin{proposition}
The boundary homeomorphism $\varphi: \partial \dd \to \partial \dd$ satisfies the $p$-Douglas condition for $p > 1 $ if and only if
\[\sum_{n = 1}^\infty  \sum_{k=1}^{2^n} |\varphi(I_{N,m})|^p 2^{(p-2)N}<\infty.\]
\end{proposition}

For the proof, we refer to \cite{KKO}.

\section{The bi-Sobolev extension}

In this section we prove Theorem \ref{thm:main}.

\begin{proof}
It is clear that the 'if' part of the statement is a trivial consequence of the fact that the $p$-Douglas condition \eqref{eq:pdouglas} characterizes the trace space of $\W^{1,p}(\dd)$. Hence we prove here the 'only if' part by assuming that $\varphi$ satisfies the $p$-Douglas condition, $\varphi^{-1}$ satisfies the $q$-Douglas condition, and constructing the desired extension $h$.

For the purposes of easier presentation we consider the boundary of the unit disk as locally flat, and hence we suppose that $\varphi$ is defined as a map of the interval $I := [0,1]$ on the real line to itself. The extension $h$ will be constructed to the upper half plane above this interval. Generally speaking, once this local extension process has been defined the map may be defined globally by a simple gluing process, see e.g. \cite{KKO} for an example of this argument.


For $n = 1,\ldots$ we let $I^{(n)} := I \times \{2^{-n}\}$ denote the unit length line segment obtained by lifting the interval $I$ to height $2^{-n}$. These segments will be defined both on the domain and target side, but to avoid confusion we change the notation to use $J$ instead of $I$ when we are on the target side. For each $n$ we also lift the boundary map $\varphi$ to a map $\varphi^{(n)}$ which maps $I^{(n)}$ to $J^{(n)}$.

The map $h$ will be constructed in a way that it maps each segment $I^{(n)}$ to the corresponding segment $J^{(n)}$. On each of these segments $h$ will be a piecewise linear approximation of $\varphi^{(n)}$, and in the strip between $I^{(n)}$ and $I^{(n+1)}$ the map $h$ will be defined by a simple gluing process we define later.

\begin{figure}
\includegraphics[scale=0.4]{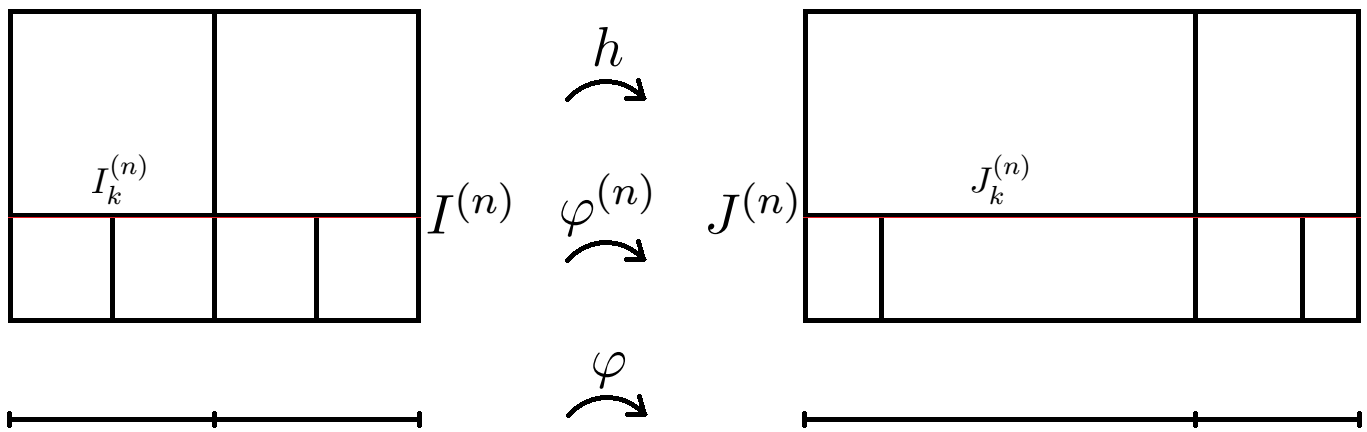}
\end{figure}

Let us start by decomposing $I^{(n)}$ dyadically into segments $I^{(n)}_k, k = 1, \ldots, 2^n$ of length $2^{-n}$, ordered from left to right. The corresponding segments $J^{(n)}_k \subset J^{(n)}$ on the target side are defined as the images of $I^{(n)}_k$ under $\varphi^{(n)}$, or in other words by projecting each $I^{(n)}_k$ down to the interval $I$, applying the boundary map $\varphi$, and then lifting the image interval back to the height $2^{-n}$.

Now it would be simple if we could just define $h$ to map each $I^{(n)}_k$ to each $J^{(n)}_k$ linearly, but it turns out that this choice is only compatible with Sobolev-estimates for $h$ and not $h^{-1}$. The issue is when the image segment $J^{(n)}_k$ happens to be quite small which results in large stretching for $h^{-1}$. To fix this issue, we define the following process.

We define some number of new segments $S^{(n)}_j$, each of which will be a union of one or more successive segments $I^{(n)}_k$ (plus one half interval in some exceptional cases). The segments $S^{(n)}_j$ are defined inductively as follows. We start by choosing $k_1$ as the smallest positive integer so that the union $J^{(n)}_1 \cup \cdots \cup J^{(n)}_{k_1}$ has length at least $2^{-n}$, and define $S^{(n)}_1$ as the union
\[S^{(n)}_1 = I^{(n)}_1 \cup \cdots \cup I^{(n)}_{k_1}.\]
We also define $T^{(n)}_1 =  J^{(n)}_1 \cup \cdots \cup J^{(n)}_{k_1}$.

Then we continue this process by defining $k_{j}$ as the smallest positive integer so that the union $J^{(n)}_{k_{j-1} + 1} \cup \cdots \cup J^{(n)}_{k_{j}}$ has length at least $2^{-n}$, and define
\[S^{(n)}_{j} = I^{(n)}_{k_{j-1}+1} \cup \cdots \cup I^{(n)}_{k_{j}} \quad \text{ and } \quad T^{(n)}_{j} = J^{(n)}_{k_{j-1}+1} \cup \cdots \cup J^{(n)}_{k_{j}}.\]
However, we add two exceptions to this. If $k_j > 1$ and the last interval $J^{(n)}_{k_j}$ to be added has length larger than $2 \cdot 2^{-n}$, then instead of adding the whole segment $I^{(n)}_{k_j}$ to $S^{(n)}_j$  we split $I^{(n)}_{k_j}$ into two halves $I^{(n)}_{k_j,-}$ and $I^{(n)}_{k_j,+}$ from left to right. We also split $J^{(n)}_{k_{j}}$ into two intervals $J^{(n)}_{k_j,-}$ and $J^{(n)}_{k_j,+}$ from left to right, where the split is chosen so that $|J^{(n)}_{k_j,-}| = 2^{-n}$. Finally we define
\[S^{(n)}_{j} = I^{(n)}_{k_{j-1}+1} \cup \cdots \cup I^{(n)}_{k_j,-}\quad \text{ and } \quad T^{(n)}_{j} = J^{(n)}_{k_{j-1}+1} \cup \cdots \cup J^{(n)}_{k_j,-}\]
instead. Furthermore, in this situation we also define the next pair of segments with index $j+1$ by $S^{(n)}_{j+1} := I^{(n)}_{k_j,+}$ and $T^{(n)}_{j+1} := J^{(n)}_{k_j,+}$. This exception guarantees that we do not combine too large intervals with small ones while keeping the length of both $S^{(n)}_{j}$ and $T^{(n)}_{j}$ larger than $2^{-n}/2$ for all $j$.

\begin{figure}
\includegraphics[scale=0.6]{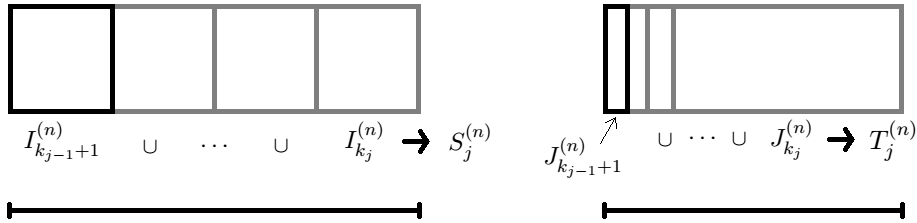}
\end{figure}

The second exception we place concerns the case of trying to define the final segment $S^{(n)}_j$, as in this case it might not be possible to choose a number $k_j$ so that the union $J^{(n)}_{k_{j-1} + 1} \cup \cdots \cup J^{(n)}_{k_{j}}$ has length at least $2^{-n}$ if there is not enough space left in $J^{(n)}$. In this case we simply forget about this segment and add the remaining segments $I^{(n)}_{k_{j-1}+1}, \ldots, I^{(n)}_{2^n}$ to the end of the previous segment $S^{(n)}_{j-1}$ instead. On the target side we similarly add $J^{(n)}_{k_{j-1} + 1} \cup \cdots \cup J^{(n)}_{k_{j}}$ to the previous segment $T^{(n)}_{j-1}$ in this case.

This process guarantees that $I^{(n)}$ has been decomposed into a number of segments $S^{(n)}_j$ of length at least $2^{-n}/2$, each of whose ''image'' segment $T^{(n)}_j$ has length at least $2^{-n}$. The map $h$ is now finally defined on $I^{(n)}$ so that it maps each $S^{(n)}_j$ to $T^{(n)}_j$ linearly.

We then extend $h$ to the strip between $I^{(n)}$ and $I^{(n+1)}$ for each $n$. Let $I^*$ be one of the segments $I^{(n)}_k$ in $I^{(n)}$ and let $I^*_-$ and $I^*_+$ denote the two intervals in $I^{(n+1)}$ which have half the length of $I^*$ and lie directly below it (the ''dyadic children'' of $I^*$, at least in projection). Let $X_1,X_2$ denote the left and right endpoint of $I^*$ and $Y_1,Y_2,Y_3$ denote the endpoints of $I^*_-$ and $I^*_+$ from left to right with $Y_2$ being the common endpoint in-between.

We now triangulate the quadrilateral $X_1 Y_1 Y_3 X_2$ via the segments $X_1 Y_1$, $X_1 Y_2$, $X_2 Y_2$ and $X_2 Y_3$ into three triangles $\Delta_i, i = 1,2,3$, on each of which we define the map $h$ as a linear map into the corresponding image triangle $\Delta_i'$. This image triangle is already well-defined since $h$ is defined on all of the vertices $X_i, i = 1,2$ and $Y_j, j = 1,2,3$. It remains to control the behaviour of $h$ on $\Delta_i$.

The triangle $\Delta_i$ on the domain side is bilipschitz-equivalent to a right triangle with sides of length $2^{-n}$. The image triangle $\Delta_i'$ has height $2^{-n}$, has a horizontal side whose length is $h(I')$, where $I'$ denotes one of the intervals $I^*, I^*_-,$ or $I^*_+$ depending on choice of $i$. In the earlier process of joining intervals $J^{(n)}_k$ into the intervals $T^{(n)}_j$, we may have shifted the respective endpoints so that the intervals $h(I')$ and $\varphi(I')$ are not equal (abusing notation slightly as $\varphi$ should be either $\varphi^{(n)}$ or $\varphi^{(n+1)}$ here). However, we claim that $|h(x) - \varphi(x)| \leq C 2^{-n}$ for a uniform constant $C$, meaning that the endpoints are shifted by an amount comparable to $2^{-n}$.

This is due to the fact that we only join multiple intervals $J^{(n)}_k$ together within the interval $T^{(n)}_j$ which already has length comparable to $2^{-n}$. And in some special cases where  $J^{(n)}_k$ is long, either the endpoints are left alone or shifted by exactly $2^{-n}$. Thus our claim is true, and this allows us compose with a bilipschitz linear shear map to assume that the image triangle $\Delta_i'$ is actually a right triangle with one side $h(I')$ and height $2^{-n}$.

Hence we need to estimate the differential and distortion of a linear map between two right triangles: one with sides of length $2^{-n}$ and another with side lengths $2^{-n}$ and $h(I')$. In brief terms,
\begin{align*}&\int_{\Delta_i} |Dh|^p \, dz \approx |h(I')|^p 2^{-n(2-p)} + 2^{-2n} \qquad \text{ and } \\& \int_{\Delta_i'} |Dh^{-1}|^q \, dz \approx |h(I')| 2^{-n}\frac{2^{-qn}}{|h(I')|^{q}} + 2^{-2n}.
\end{align*}
Note that the expression $2^{-2n}$ is summable over all dyadic intervals, which means that we need only control the first terms in each expression. Let us begin by bounding the $p$-norm of $Dh$.

The term $|h(I')|^p 2^{-n(2-p)}$ is only relevant when $|h(I')|$ is much larger than $2^{-n}$. In this case the interval $h(I')$ was obtained from $\varphi(I')$ in the construction either by changing nothing or removing a small segment. In either case $|h(I')| \leq |\varphi(I')|$. But now it only remains to control the sum
\[\sum_{n=1}^{\infty} \sum_{k=1}^{2^n} |\varphi(I_{n,k})|^p 2^{-n(2-p)},\]
which is known to be finite since $\varphi$ satisfies the $p$-Douglas condition and thus also the discrete $p$-Douglas condition.

Let us then bound the $q$-norm of $Dh^{-1}$. It is enough to bound the first term $2^{-n(1+q)}|h(I')|^{1-q}$. Let us suppose that the segment $I'$ is part of the segment $S^{(n)}_j$, and that $S^{(n)}_j$ consists of the union of $N \geq 1$ neighbouring segments of the same length. Due to the special cases in the construction above, it may be possible that one of these segments is half the length of the others, but this will not particularly affect the estimates so we disregard this case. We estimate the total energy coming from all these segments, which equals to
\begin{equation}\label{eq:qstuff}
N \cdot 2^{-n(1+q)}|h(I')|^{1-q} = 2^{-n}\frac{(N2^{-n})^{q}}{(N|h(I')|)^{q-1}},
\end{equation}
as the map $h$ is linear on $S^{(n)}_j$ and thus maps each of the $N$ segments into a segment of the same length $|h(I')|$. The image segment $h(S^{(n)}_j) = T^{(n)}_j$ was chosen crucially to have length larger than $2^{-n}$ but less than $2^{-n+1}$. Due to this fact we may argue as follows:

Each segment $T^{(n)}_j$ can be covered by at most three neighbouring dyadic segments $J^{(n)}_l$ on the image side. For each of these segments on the image side, we recall that sum in the   discrete $q$-Douglas condition for $\varphi^{-1}$ involves a corresponding term of the form
\[2^{(q-2)n} |\varphi^{-1}(J^{(n)}_l)|^q.\]
Notice that $2^{-n} \leq |T^{(n)}_j| = N|h(I')|$. Moreover, since $|I'| = 2^{-n}$ it holds that $N 2^{-n} = |S^{(n)}_j|$. But since $T^{(n)}_j$ was covered by the intervals $J^{(n)}_l$, of which there were at most three, we also have $|S^{(n)}_j| \leq 3 \max_{l} |\varphi^{-1}(J^{(n)}_l)|$. Thus, estimating the energy from \eqref{eq:qstuff}, we get:
\[2^{-n}\frac{(N2^{-n})^{q}}{(N|h(I')|)^{q-1}} \leq C 2^{-n} \frac{\max_{l} |\varphi^{-1}(J^{(n)}_l)|^q}{(2^{-n})^{q-1}}.\]
Considering also that each dyadic interval $J^{(n)}_l$ is involved in this process at most three times, we may sum over $n$ and $l$ to obtain that
\[\int_{\dd} |Dh^{-1}(z)|^q \, dz \leq C \sum_{n=1}^{\infty} \sum_{l = 1}^{2^n} 2^{(q-2)n} |\varphi^{-1}(J^{(n)}_l)|^q.\]
But this is finite due to the discrete $q$-Douglas condition for $\varphi^{-1}$. This finishes the proof.

\end{proof}

\section{Sobolev-estimates for the harmonic extension} \label{sec:harm}

In this section we give the proof of Theorem \ref{thm:rkcp>2}, showing that the harmonic extension will lie in $\W^{1,p}(\dd)$ as long as the boundary map $\varphi : \partial \dd\to\partial \dd$ satisfies the $p$-Douglas condition. We follow the same line of arguments as in \cite{HKW}, but provide the details here for the readers' convenience.

Split the unit circle into dyadic arcs $I_{n,k}$. For each given arc $I_{n,k}$, we now inductively define a family $\P(I_{n,k})$ of disjoint dyadic arcs that cover the whole circle as follows. We first add $I_{n,k}$ into $\P(I_{n,k})$. We then repeat the following process:

Recall that every dyadic arc $J$ has a unique sibling $J'$ being the unique arc for which $J \cup J'$ is also a dyadic arc. For each arc $J$ in $\P(I_{n,k})$, we now consider two neighbouring arcs of $J$. The first of these is the unique sibling $J'$, which we add to $\P(I_{n,k})$ if $J'$ was not already contained in the union of all arcs added to $\P(I_{n,k})$ so far. We also let $J^*$ denote the dyadic arc which neighbours $J$ and has twice the size of $J$, which necessarily lies on the other side as $J'$. We also add $J^*$ to $\P(I_{n,k})$ unless $J*$ was already contained in the union of previous arcs.

The key properties of the family $\P(I_{n,k})$ is that the arcs in it are disjoint, cover the whole boundary, and the distance from each arc $J \in \P(I_{n,k}), J \neq I_{n,k}$ to the midpoint of $I_{n,k}$ is comparable to the length of $J$.

We now estimate the differential of the harmonic extension. The Poisson extension formula may be differentiated to obtain that
\[h_z(z) = \frac{1}{2\pi} \int_0^{2\pi} \frac{e^{i\theta}}{(z-e^{i\theta})^2} \varphi(e^{i\theta}) \, d\theta.\]
Let us write $\varphi(e^{i\theta}) = e^{i f(\theta)}$ for a real-valued continuous increasing function $f$. We then integrate by parts to estimate as follows
\begin{align*}
|h_z(z)| &= \left| \frac{1}{2\pi} \int_0^{2\pi} \frac{e^{i\theta}}{(z-e^{i\theta})^2} e^{i f(\theta)} \, d\theta \right|
\\&= \left| \frac{1}{2\pi i} \int_0^{2\pi} \frac{1}{(z-e^{i\theta})} \, d\left(e^{i f(\theta)}\right) \right|
\\
&\leq \frac{1}{2\pi} \int_0^{2\pi} \frac{1}{|z-e^{i\theta}|} d \mu_f(\theta),
\end{align*}
where $\mu_f$ denotes the Lebesgue-Stieltjes measure of $f$.

Let us estimate this expression for $z \in U_{n,k}$, using the family $\P(I_{n,k})$ to separate the domain of integration so that we may estimate this expression separately for each $e^{i\theta} \in J \in \P(I_{n,k})$. On each such arc $J$, the distance $|z - e^{i\theta}|$ is comparable to the length of $J$, which is equal to $2^{-N}$ for some $N$ with $N \leq n$. Thus
\[
\left|\int_{\{0\leq\theta\leq 2\pi, e^{i\theta}\in J\}} \frac{1}{|z-e^{i\theta}|} d \mu_f(\theta)\right|
\leq \frac{|\varphi(J)|}{2^{-N}}.
\]
Let us then write $J = I_{N,m}$ for some $m$, and note that $N \leq n$. Thus for $z \in U_{n,k}$ we have
\[|h_z(z)| \leq \sum_{N \leq n}\ \sum_{I_{N,m}\in\P(I_{n,k})} |\varphi(I_{N,m})| 2^N.\]
We raise to the power $p$, integrate over $U_{n,k}$, and sum over $n$ and $k$ to obtain
\begin{align*}
\int_{\dd} |h_z(z)|^p \, dz &\leq
\sum_{n=1}^{\infty} \sum_{k=1}^{2^n} 2^{-2n} \left( \sum_{N \leq n}\ \sum_{I_{N,m}\in\P(I_{n,k})} |\varphi(I_{N,m})| 2^N \right)^p
\\&\leq
\sum_{n=1}^{\infty} \sum_{k=1}^{2^n} 2^{-2n} \left(\sum_{N \leq n}\ \sum_{I_{N,m}\in\P(I_{n,k})} \frac{1}{\left(2^{-\alpha N}\right)^{\frac{p}{p-1}}}\right)^{p-1}
\\& \qquad \qquad \cdot \sum_{N \leq n}\ \sum_{I_{N,m}\in\P(I_{n,k})} |\varphi(I_{N,m})|^p 2^{(1-\alpha)pN}
\end{align*}
Here we have applied H\"older's inequality in the second step for the sequences $2^{\alpha N}$ and $|\varphi(I_{N,m})|2^{(1-\alpha)N}$, where $\alpha \in (0,1/p)$ is a parameter which will matter later. For each $N \leq n$, the number of $m$ for which $I_{N,m} \in \P(I_{n,k})$ is definitely less than $5$. Thus in the first factor above we may estimate the sum as a geometric series where only the leading term matters:
\[\left(\sum_{N \leq n}\sum_{I_{N,m}\in\P(I_{n,k})} \frac{1}{\left(2^{-\alpha N}\right)^{\frac{p}{p-1}}}\right)^{p-1} \leq C \left( \frac{1}{\left(2^{-\alpha n}\right)^{\frac{p}{p-1}}} \right)^{p-1} = C 2^{\alpha pn}.\]
Let $\chi_{N,m,n,k} := \chi\left(I_{N,m}\in\P(I_{n,k})\right)$ denote a function which is $1$ when $I_{N,m}\in\P(I_{n,k})$ and $0$ otherwise. Then, combining with the above estimates, we obtain

\begin{align*}
\int_{\dd} |h_z(z)|^p \, dz &\leq
C \sum_{n=1}^{\infty} \sum_{k=1}^{2^n} 2^{\left(\alpha p-2\right)n} \sum_{N \leq n}\sum_{m=1}^\infty |\varphi(I_{N,m})|^p 2^{(1-\alpha)pN} \chi_{N,m,n,k}
\\&=
C \sum_{n=1}^{\infty}\sum_{N \leq n}\sum_{m=1}^\infty  2^{\left(\alpha p-2\right)n} |\varphi(I_{N,m})|^p 2^{(1-\alpha)pN} \sum_{k=1}^{2^n} \chi_{N,m,n,k}
\\&=
C \sum_{N = 1}^\infty \sum_{n=N}^{\infty} \sum_{m=1}^\infty  2^{\left(\alpha p-2\right)n} |\varphi(I_{N,m})|^p 2^{(1-\alpha)pN} \sum_{k=1}^{2^n} \chi_{N,m,n,k}
\end{align*}
We should now estimate the sum $\sum_{k=1}^{2^n} \chi_{N,m,n,k}$ as a function of $N,m$, and this sum represents the amount of intervals $I_{n,k}$, $k = 1,\ldots,2^n$ for which $I_{N,m}\in\P(I_{n,k})$. As before, for a fixed $k$ the number of $m$ for which this relation holds is bounded by $5$. But conversely, as the amount of $k$ for fixed $m$ must be the same for each $m$ by symmetry, we find that the exact amount of such $k$ must be comparable to $2^{n-N}$. Thus
\begin{align*}
\int_{\dd} |h_z(z)|^p \, dz &\leq
C \sum_{N = 1}^\infty \sum_{n=N}^{\infty} \sum_{m=1}^\infty  2^{\left(\alpha p-2\right)n} |\varphi(I_{N,m})|^p 2^{(1-\alpha)pN} 2^{n-N}
\\&=
C \sum_{N = 1}^\infty  \sum_{m=1}^\infty |\varphi(I_{N,m})|^p 2^{(1-\alpha)pN-N} \sum_{n=N}^{\infty} 2^{\left(\alpha p-1\right)n}
\\&\leq
C \sum_{N = 1}^\infty  \sum_{m=1}^\infty |\varphi(I_{N,m})|^p 2^{(1-\alpha)pN-N} 2^{\left(\alpha p-1\right)N}
\\&=
C \sum_{N = 1}^\infty  \sum_{m=1}^\infty |\varphi(I_{N,m})|^p 2^{(p-2)N},
\end{align*}
where we used the fact that $\alpha$ was chosen in a way that $\alpha p - 1 < 0$ to calculate the geometric series. The finiteness of the final sum is exactly the discrete $p$-Douglas condition of $\varphi$, which proves the claim.

\section{Failure of the harmonic extension in the bi-Sobolev case} \label{sec:harmbad}

In this section we prove the statement of Example \ref{ex:biharm}, showing that the harmonic extension does not solve the bi-Sobolev extension problem.

\begin{proof} Let $\varphi : \partial \dd \to \partial \dd$ be a boundary map to be chosen. We wish to pick $\varphi$ so that the distortion function $K_h$ of the harmonic extension $h$ has a singularity at $z = 1$ which will result in $K_h$ not belonging to $L^1(\dd)$. The map $\varphi$ will be chosen as a smooth homeomorphism of the boundary to itself.

Due to Heinz \cite{HZ}, the differential of $h$ is strictly bounded from below in the whole disk $\dd$. Moreover, due to the smoothness of $h$ up to the boundary the differential is also uniformly bounded from above. Let $c$ and $C$ be constants so that $c \leq |Dh| \leq C$. Then for $r = |z| > 1/2$ we can estimate the distortion from below:
\[K_h = \frac{|Dh|^2}{J_h} = \frac{|Dh|^2}{\frac{1}{r} \Im m (h_r \bar{h_\theta})} \geq \frac{c^2}{2C} \frac{1}{|h_\theta|}.\]
The Poisson integral formula may be used to find an integral representation for the angular derivative:
\[h_\theta(z) = \int_{\partial \dd}\frac{1 - |z|^2}{|z-\omega|^2} \varphi_\theta(\omega) |d \omega|.\]
We now pick the boundary map $\varphi$ as the map
\[\varphi(e^{i\theta}) = \exp\left(i\pi \, e^{1-(\pi/\theta)^4}\right).\]
We consider $\varphi$ in the regions
\[S_n = \{r e^{i\theta} : 2^{-n-1} < 1-r < 2^{-n},\  |\theta| < \frac{1}{n^{0.25}}\}.\]
If we can show that $|h_\theta| \leq C n^{0.5}2^{-n}$ in $S_n$, then
\begin{equation}\label{eq:someK}\int_{S_n} K_h \, dz \geq \int_{S_n} \frac{c_1}{|h_\theta|} \, dz \geq c_2 \frac{2^{-n}}{n^{0.25}} \frac{1}{n^{0.5} 2^{-n}} = \frac{c_2}{n^{0.75}}.\end{equation}
Since the sets $S_n$ are disjoint, taking the sum over all $n$ will show that $\int_{\dd} K_h \, dz = \infty$ as desired. For $z \in S_n$ we therefore estimate that
\begin{align*}
|h_\theta(z)| &= (1 - |z|^2) \left|\int_{-\pi}^{\pi}\frac{1}{|z-e^{it}|^2} \frac{i 4 \pi^5}{t^5}e^{1-(\pi/t)^4}\exp\left(i\pi \, e^{1-(\pi/t)^2}\right) \, dt \right|
\\&\leq C_1 2^{-n} \int_{-\pi}^{\pi}\frac{1}{|z-e^{it}|^2} \frac{1}{t^5} e^{-(\pi/t)^4}\, dt
\\&\leq C_2 2^{-n} \int_0^{\pi}\frac{1}{|z-e^{it}|^2} e^{-(2/t)^4}\, dt.
\end{align*}
Let us now split the region of integration into two intervals. First, we consider $t \in [0,2/n^{0.25}]$. Here we may use the estimates $|z - e^{it}| \geq 2^{-n}$ and $(2/t)^4 \geq n$ to get
\[2^{-n} \int_0^{2/n^{0.25}}\frac{1}{|z-e^{it}|^2} e^{-(2/t)^4}\, dt
\leq \frac{1}{n^{0.25}} 2^{n} e^{-n}
\leq \frac{1}{n^{0.25}}.
\]
For $t > 2/n^{0.25}$, we note that $|z - e^{it}| \geq c/n^{0.25}$. Thus we find that
\[2^{-n} \int_{2/n^{0.25}}^\pi\frac{1}{|z-e^{it}|^2} e^{-(2/t)^4}\, dt \leq 2^{-n} \frac{1}{(c/n^{0.25})^2} = \frac{1}{c} n^{0.5}2^{-n}.\]
This proves the claim through \eqref{eq:someK}.
\end{proof}

\section{Sharpness of  Beurling-Ahlfors type extensions, Example~\ref{ex:p-integrabledist}}

In this section we first construct Example~\ref{ex:p-integrabledist}, showing that the Karafyllas-Ntalampekos condition~\eqref{eq:p-integralbdd} is not a necessary condition for a boundary map to admit a homeomorphic extension with $p$-integrable distortion.

Fix $p \ge 1$ and $q>p$. Then we choose $\epsilon \in (0, q-p)$. Consider the radial stretching mapping $h \colon \mathbb C  \to \mathbb C$ defined by
\[h(z) = H(|z|) \frac{z}{|z|}, \qquad \textnormal{ where} \quad  H(s) = \exp\left(- s^{-\frac{2}{p+\epsilon}}\right).\]
It will be enough to show that $K_h\in \mathscr L^p$ near $z=0$   but that the restriction $\varphi = h \vert_{\rr} \colon \R \onto \R$ of $h$ to the real line satisfies
\begin{equation}\label{eq:notalampekos}
\int_0^1 \int_0^1 \log^q \left(e+  \frac{|\varphi(x+t) - \varphi(x)|}{|\varphi(x) - \varphi(x-t)|}  \right) \, dx\, dt \ = \ +\infty.
\end{equation}
We simplify writing and denote  $\beta = \frac{2}{p+\epsilon}$.  A direct computation shows that 
\[\abs{Dh(z)}^2 = \max \left\{ \dot{H}^2 (\abs{x}), \frac{H^2(\abs{x})}{\abs{x}^2} \right\} \le \frac{H^2(\abs{z})}{\abs{z}^2} \left(\beta^2 \abs{z}^{-2\beta}  +1\right)   \]
and
\[J_h(z) = \dot{H} (\abs{z}) \frac{H(\abs{z})}{\abs{z}} = \beta \,  \frac{H^2(\abs{z})}{\abs{z}^{2+\beta}} \, . \]
Therefore,
\[K_h(z) \le \beta \abs{z}^{-\beta} + \beta^{-1} \abs{z}^\beta  \, . \]

Let $\mathbb B_R$ be the ball centered at $0$ with radius $R>0$. Computing in polar coordinates, we obtain that
\[\int_{\mathbb B_R} K_h^p(z) \, dz = 2\pi \int_0^R \left( \beta  r^{1 - \frac{2p}{p+\epsilon}} + \beta^{-1} r^{1+   \frac{2p}{p+\epsilon}   }    \right)\, \dtext  r < \infty.\]
Let us now show that \eqref{eq:notalampekos} holds. Consider the subset $U = \{(x,t) : t \geq 2x\}$ of $[0,1]^2$. Note that
\[\varphi'(s) = \dot{H}(s) = \beta  s^{-\beta-1}\exp\left(- s^{-\beta }\right).\]
Thus both $\varphi$ and $\varphi'$ are increasing for $s$ close to zero. For all  $(x,t) \in U$ close to the origin we then estimate using the mean value theorem that 
\begin{align}\label{eq:abovevarphi}
\frac{\varphi(x+t) - \varphi(x)}{\varphi(x) - \varphi(x-t)}
&\geq \frac{\varphi(x+t) - \varphi(x+t/2)}{\varphi(x) - \varphi(x-t)} \nonumber
\\  &\geq \frac12 \frac{\varphi'(x + t/2)}{\varphi'(x)}
\\&\geq \frac12 \frac{\varphi'(2x)}{\varphi'(x)}. \nonumber
\end{align}
We may now simplify this expression to find that
\[\frac12\frac{\varphi'(2x)}{\varphi'(x)} \geq C x^\alpha \exp\left(\left(1-2^{-\beta }\right)x^{-\beta}\right),\]
where $C, \alpha > 0$ are constants. Note that $1-2^{-\beta} > 0$.  Since $\beta=\frac{2}{p+\varepsilon} < \frac{2}{q}$ we may further estimate this expression for $x$ close to zero by using a single exponential to dominate the lower order factors:
\[\frac{\varphi'(2x)}{\varphi'(x)} \geq C x^\alpha \exp\left(\left(1-2^{-\beta}\right)x^{-\beta}\right) \geq  \exp\left(x^{-\frac{2}{q}}\right),\]
for all $0<x<\alpha_{p,q}=\alpha$. Hence
\begin{align*}
\int_0^1 \int_0^1 \log^q \left(e+  \frac{|\varphi(x+t) - \varphi(x)|}{|\varphi(x) - \varphi(x-t)|}  \right)   \, dx\, dt & \geq  c \int_0^\alpha  x^{-1} \, \dtext x = \infty \, . 
\end{align*}
This concludes the construction of Example~\ref{ex:p-integrabledist}.

Question 1.5 in~\cite{KN}  also   asks if the  boundary condition
\begin{equation}
\int_0^{2\pi }  \int_0^{2\pi } \exp\big(q \delta_\varphi (\theta, t)\big) \,  \dtext \theta \,  \dtext t   < \infty \, , \qquad \textnormal{where } q>0
\end{equation} 
is necessary for obtaining an extension $h \colon \mathbb D \onto \mathbb D$ of exponentially integrable distortion; that is, $\exp(K_h) \in L^\lambda (\mathbb D)$ for some $\lambda >0$. To see that this is not the case, for fixed $\lambda >0$ we may choose
\[H(s) = \exp \left(- \mu \log^2 (e/s) \right) \qquad \textnormal{ where } 0< 2 \mu < \lambda \, ,    \]
and consider the corresponding radial symmetric map $h= H(\abs{z}) \frac{z}{\abs{z}}$.
Since
\[K_h (z) \le 2 \mu \log (e/\abs{z}) + 2\mu^{-1} \log^{-1} (e/\abs{z})   \]
we have that $\exp( \lambda K_h) \in L_{\loc}^1 (\mathbb C)$. On the other hand,  the restriction $\varphi = h \vert_{\rr} \colon \R \onto \R$ of $h$ to the real line satisfies
\begin{equation}\label{eq:notalampekos2}
\int_0^1 \int_0^1 \exp \left( q \cdot  \frac{|\varphi(x+t) - \varphi(x)|}{|\varphi(x) - \varphi(x-t)|}  \right) \, \dtext x\, \dtext t \ = \ +\infty \qquad \textnormal{for all } q>0 \, . 
\end{equation}
Indeed, estimating as in~\eqref{eq:abovevarphi}, for all $(x,t) \in U$ close to the origin, we have
\[
\begin{split}
\frac{\varphi(x+t) - \varphi(x)}{\varphi(x) - \varphi(x-t)}
&\geq \frac{1}{2} \frac{\log (1/(2x))H({2x})  }{\log (1/(x))H({x}) } \ge \frac{1}{\log (1/x)} \frac{H(2x)}{H(x)} \\
& \geq \frac{1}{\log (1/x)} \frac{1}{(2x)^\mu}
\end{split}
\]
and therefore ~\eqref{eq:notalampekos2} follows.

\end{document}